\documentclass[12pt]{article}
\usepackage{amssymb}

\usepackage{graphicx}
\usepackage{amsmath}

\advance\voffset-1truecm\relax
\advance\hoffset-1truecm\relax
\font\chu=cmr10

\begin{document}

\title{AN EXTENSION OF\\
 UNIQUENESS THEOREMS\\
 FOR MEROMORPHIC MAPPINGS }
\author{$\quad $Gerd Dethloff and Tran Van Tan}
\date{$\quad$}
\maketitle

\noindent 
\begin{abstract}  %
%

\noindent  \textit{In this paper, we give some results on the number of 
meromorphic
mappings of} $\mathbb{C}^{m}$ \textit{ into} $\mathbb{C}P^{n}$ 
\textit{under a condition on the
inverse images of hyperplanes in} $\mathbb{C}P^{n}$. \textit{At the 
same time, we
give an answer for an open question posed by H. Fujimoto in 1998.}

\end{abstract}  %
%

\section{\textbf{\ Introduction } \vskip0.2cm}

In 1926, R. Nevanlinna showed that for two nonconstant meromorphic
functions $f$ and g on the complex plane $\mathbb{C}$, if they have the 
same
inverse images for five distinct values, then $f=g$, and that $g$ is a 
special
type of a linear fractional tranformation of $f$ if they have the same
inverse images, counted with multiplicities, for four distinct values.

In 1975, H. Fujimoto [2] generalized Nevanlinna's result to the case of 
meromorphic mappings of $\mathbb{C}^{m}$ into $\mathbb{C}P^{n}.$ This 
problem continued to be studied by L. Smiley $[9],$ S.Ji $ [5]$ and 
others.

Let $f$ be a meromorphic mapping of $\mathbb{C}^{m}$ into 
$\mathbb{C}P^{n}$
and $H$ be a hyperplane in $\mathbb{C}P^{n}$ such that im$f\nsubseteq 
H.$
Denote by $v_{(f,H)}$ the map of $\mathbb{C}^{m}$ into $\mathbb{N}_{0}$ 
such
that $v_{(f,H)}(a)$ $(a\in \mathbb{C}^{m})$ is the intersection 
multiplicity
of the image of $f$ and $H$ at $f(a).$ Let $k$ be a positive interger 
or $%
+\infty .$ We set

\hbox to 5cm {\hrulefill }

\noindent
{\chu 2000 Mathematics Subject Classification:} 
{\chu  32\,H 30.}

\noindent
{\chu Key words and phrases:} {\chu uniqueness theorem, meromorphic 
mapping, linearly degenerate.}

\newpage

\begin{equation*}
v_{(f,H)}^{k)}(a)=\begin{cases} 0\;\;\;\;\;\;\;\;\;\;\;\text{ 
if}\;\;\;v_{(f,H)}(a)>k,\\
v_{(f,H)}(a)\text{ if}\;\;\;v_{(f,H)}(a)\le k.
\end{cases}
\end{equation*}

Let $f$ be a linearly nondegenerate meromorphic mapping of 
$\mathbb{C}^{m}$
into $\mathbb{C}P^{n}$ and $\{H_{j}\}_{j=1}^{q}$ be $q$ hyperplanes in
general position with

\noindent (a) \;\;\; \textit{dim} $\Big\{z:v_{(f,H_i)}^{k)}(z)>0$ 
\;\;and\;\; $%
v_{(f,H_j)}^{k)}(z)> 0\Big\}\leq m-2$
\noindent for all $1\leq i < j \leq q.$

For each positive integer $p$, denote by 
$F_{k}(\{H_{j}\}_{j=1}^{q},f,p$ $)$
the set of all linearly nondegenerate meromorphic mappings $g$ of 
$\mathbb{C}%
^{m}$ into $\mathbb{C}P^{n}$ such that:

\noindent (b)\;\;\;\;\; \textit{min}$\big\{ v_{(g, H_j)}^{k)}, p\big\} 
=\textit{min}\big\{ %
v_{(f, H_j)}^{k)}, p\big\},$ 

\noindent (c) \ \ \ \ \ $g=f$ on $\bigcup\limits_{j=1}^{q}\big\{%
z:v_{(f,H_{j})}^{k)}(z)>0\big\}.$

In $[5]$, S.Ji showed the following\\

\noindent%
%
\textbf{Theorem J.} ([5])
\textit{If }$q=3n+1$ \textit{\ and }$k={+\infty },$\textit{\ then for 
three mappings }$%
f_{1},f_{2},f_{3}\in $\textit{\ 
}$F_{k}\big(\{H_{j}\}_{j=1}^{q},f,1\big)
,$\textit{\ the mapping }$f_{1}\times f_{2}\times f_{3}:\mathbb{C}%
^{m}\longrightarrow \mathbb{C}P^{n}\times \mathbb{C}P^{n}\times 
\mathbb{C}%
P^{n}$\textit{\ is algebraically degenerate, namely, }$%
\{(f_{1}(z),f_{2}(z),f_{3}(z))$\textit{, }$z\in 
\mathbb{C}^{m}\}$\textit{\
is contained in a proper algebraic subset of }$\mathbb{C}P^{n}\times 
\mathbb{C%
}P^{n}\times \mathbb{C}P^{n}.$\\

In 1929, H. Cartan declared that there are at most two meromorphic 
functions
on $\mathbb{C}$ which have the same inverse images (ignoring 
multiplicities)
for four distinct values. However in 1988, N. Steinmetz ([10]) gave 
examples which
showed that H. Cartan's declaration is false. On the other hand, in 
1998, Fujimoto ([4])
showed that H. Cartan's declaration is true if we assume that
meromorphic functions on $\mathbb{C}$ share four distinct values 
counted
with multiplicities truncated by 2. He gave the following theorem\\

\noindent%
%
\textbf{Theorem F.} ([4])
\textbf{\ }\textit{If }$q=3n+1$\textit{\ and }$k={+\infty }$\textit{\ 
then }$F_{k}\big(%
\{H_{j}\}_{j=1}^{q},f,2\big)
$\textit{\ contains at most two mappings.}\\

He also proposed an open problem  asking if the number $q=3n+1$ in 
Theorem F can be
replaced by a smaller one. Inspired by this question, in this paper we 
will
generalize the above results to the case where the number $q=3n+1$ is 
in fact
replaced by a smaller one. We also obtain an improvement concerning
truncating multiplicities.

Denote by $\Psi $ the Segre embedding of \textit{\ 
}$\mathbb{C}P^{n}\times 
\mathbb{C}P^{n}$ into $\mathbb{C}P^{n^{2}+2n}$ which is defined by 
sending the
ordered pair $\left( (w_{0},...,w_{n}),(v_{0},...,v_{n})\right) $ to $%
(...,w_{i}v_{j},...)$ (in lexicographic order).

Let $h:\mathbb{C}^{m}\longrightarrow $\textit{\ } 
$\mathbb{C}P^{n}\times 
\mathbb{C}P^{n}$ be a meromorphic mapping. Let 
$(h_{0}:...:h_{n^{2}+2n})$
be a representation of $\Psi  \circ h$ . We say that $h$ is linearly 
degenerate (with the 
algebraic structure in \textit{\ } $\mathbb{C}P^{n}\times \mathbb{C}%
P^{n}$ given by the Segre embedding) if $h_{0},...,h_{n^{2}+2n}$ are
linearly dependent over $\mathbb{C}$ .

Our main results are stated as follows:\\

\noindent%
%
\textbf{Theorem 1.}
\textit{There are at most two distinct mappings in }$F_{k}\big(%
\{H_{j}\}_{j=1}^{q},f,p\big)$\textit{\ in each of the following cases}: 

\noindent i)\; \;\; $1\leq n\leq 3,q=3n+1,p=2$\textit{\ and }$23n\leq
k\leq {+\infty }$\textit{\ }

\noindent ii)\; \;$4\leq n\leq 6,q=3n,p=2$\textit{\ and }$\frac{%
(6n-1)n}{n-3}\leq k\leq +\infty $

 \noindent iii) $n\geq 7,q=3n-1,p=1$\textit{\ and }$\frac{%
(6n-4)n}{n-6}\leq k\leq {+\infty }$\textit{\ .\vskip0.28cm}
\noindent%
%
\textbf{Theorem 2.}
\textit{\ Assume that }$q=\left[ \frac{5(n+1)}{2}\right] 
,(65n+171)n\leq k\leq
+\infty $\textit{\ , where }$\left[ x\right] :=\max \{d\in 
\mathbb{N}:$%
\textit{\ }$d\leq x\}$ \textit{for\ a positive constant }$x. $ 
\textit{Then
one of following assertions holds :}

\noindent i)\;\; $\#F_{k}\big(\{H_{j}\}_{j=1}^{q},f,1\big)\leq 2.$

\noindent ii) \textit{ For any }$f_{1},f_{2}\in 
F_{k}\big(\{H_{j}\}_{j=1}^{q},f,1\big),$%
\textit{\ the mapping }$f_{1}\times f_{2}:\mathbb{C}^{m}\longrightarrow 
\mathbb{C}P^{n}\times \mathbb{C}P^{n}$ \textit{is linearly degenerate 
(with
the algebraic structure in }$\mathbb{C}P^{n}\times 
\mathbb{C}P^{n}$\textit{\
given by the Segre embedding).}\\

We finally remark that we obtained similar uniqueness theorems 
with moving targets in [11], but only with a bigger number of targets 
and with much bigger truncations.\\
\noindent \textbf{Acknowledgements: }The second author would like to 
thank Professor
Do Duc Thai for valuable discussions, the Universit\'{e} de Bretagne 
Occidentale for its hospitality and support, and the PICS-CNRS 
ForMathVietnam for its support.

\section{\noindent \textbf{\ Preliminaries}\vskip0.2cm}
We set $\Vert z\Vert:= (\vert z_1\vert^2 +\cdots + \vert 
z_m\vert^2)^{1/2}$
for $z=(z_1,\dots,z_m)\in \mathbb{C}^m,$ $B(r) := \big\{ z: \Vert 
z\Vert < r\big\},\;\; S(r) : =\big\{z: \big\Vert %
z\Vert =r\big\},\;\; d^c :=\displaystyle\frac{\sqrt{-1}} {4\pi} 
(\overline
\partial -\partial),$ $\upsilon:= (dd^c \Vert z\Vert^2)^{m-1}$ \;\; and 
\;\;$\sigma:= d^c \log
\Vert z\Vert^2 \land (dd^c \log\, \Vert z\Vert^2)^{m-1}$.

Let $F$ be a nonzero holomorphic function on $\mathbb{C}^{m}.$ For an 
$m$-tuple $%
\alpha :=(\alpha _{1},\dots ,\alpha _{m})$ of nonnegative integers, set 
$%
|\alpha |:=\alpha _{1}+\cdots +\alpha _{m}$ and $D^{\alpha 
}F:=\displaystyle%
\frac{\partial ^{|\alpha |}F}{\partial z_{1}^{\alpha _{1}}\ldots 
\partial
z_{m}^{\alpha _{m}}}$.
We define the map $v_F:\mathbb{C}^m\rightarrow \mathbb{N}_0$ by 
$v_F(z):=\max%
\big\{p:D^\alpha F(z) = 0$ for all $\alpha$ with $|\alpha| < p \big\}$. 
Let $%
k$ be a positive integer or $+\infty.$ Define the map $v_F^{k)}$ of 
$\mathbb{%
C}^m$ into $\mathbb{N}_0$ by \vskip-0.15cm
\begin{equation*}
v_{F}^{k)}(z):=\begin{cases}0\;\;\;\;\;\;\;\;\text{ 
if}\;\;\;v_{F}(z)>k,\\
v_{F}(z)\text{ \ if}\;\;\;v_{F}(z)\leq k.%
\end{cases}
\end{equation*}
Let $\varphi $ be a nonzero meromorphic function on $\mathbb{C}^{m}$. 
We
define the map $v_{\varphi }^{k)}$ as follows: For each $z\in 
\mathbb{C}%
^{m}$, choose nonzero holomorphic functions $F$ and $G$ on a
neighbourhood $U$ of $z$ such that $\varphi =\frac{F}{G}$ on $U$ and 
\textit{dim}$ %
\big(F^{-1}(0)\cap G^{-1}(0)\big)\leq m-2$. Then put $v_{\varphi
}^{k)}(z):=v_{F}^{k)}(z)$. Set\vskip-0.15cm\vskip-0.15cm 
\begin{equation*}
\big\vert v_\varphi^{k)} \big\vert := \overline {\big\{z: 
v_\varphi^{k)} (z)
0\big\}}. 
\end{equation*}
Define\vskip-0.15cm\vskip-0.15cm
\begin{equation*}
N^{k)} (r, v_\varphi) : =\int\limits_{1}^r \frac {n^{k)} (t)} { 
t^{2m-1}}
dt,\;\;\; (1 < r < +\infty) 
\end{equation*}
where \vskip-0.15cm\vskip-0.15cm 
\begin{equation*}
n^{k)} (t): = \int\limits_{\big\vert v_\varphi^{k)} \big\vert \cap 
B(t)}
v_\varphi^{k)} \upsilon\;\;\;\text{for}\;\;\; m \ge 2 
\end{equation*}
and\vskip-0.15cm\vskip-0.15cm
\begin{equation*}
n^{k)}(t):=\sum_{|z|\le t}v_{\varphi 
}^{k)}(z)\;\;\;\text{for}\;\;\;m=1. 
\end{equation*}
Set $N(r,v_{\varphi }):=N^{+\infty )}(r,v_{\varphi })$. For $l$  a
positive integer or $+ \infty$, set \vskip-0.15cm 
\begin{equation*}
N_{l}^{k)}(r,v_{\varphi 
}):=\int\limits_{1}^{r}\frac{n_{l}^{k)}(t)}{t^{2m-1}}%
dt,\;\;\;(1<r<+\infty ) 
\end{equation*}
where \ \ $n_{l}^{k)}(t):=\int\limits_{\big\vert v_{\varphi 
}^{k)}\big\vert %
\cap B(t)}\textit{min}\big\{ v_{\varphi }^{k)},l\big\}\upsilon 
\;\;\;\text{for}%
\;\;\;m\ge 2$ and $n_{l}^{k)}(t):=\sum_{|z|\le t}\textit{min}\big\{ 
v_{\varphi
}^{k)}(z),l\big\}\;\;\;\text{for}\;\;\;m=1.$
Set $\overline{N}(r,v_{\varphi }):=N_{1}^{+\infty )}(r,v_{\varphi })$ 
and $´\;\;%
\overline{N}^{k)}(r,v_{\varphi }):=N_{1}^{k)}(r,v_{\varphi }).$
For a closed subset $A$ of a purely $(m-1)$-dimensional analytic subset 
of $%
\mathbb{C}^{m}$ , we define
\begin{equation*}
\overset{\_\_}{N}(r,A):=\int\limits_{1}^{r}\frac{\overset{\_}{n}(t)}{t^{2m-1}%
}dt,\;\;\;(1<r<+\infty ) 
\end{equation*}
where
\begin{equation*}
\overset{\_}{n}(t):=\begin{cases}\int\limits_{A\cap B(t)}\upsilon 
\;\;\;\;\;\;\;\;\; \text{  for }m\geq 2\\
\#\left( A\cap B(t)\right) \text{ for }m=1.\end{cases} 
\end{equation*}

Let $f:\mathbb{C}^{m}\rightarrow \mathbb{C}P^{n}$ be a meromorphic 
mapping.
For arbitrarily fixed homogeneous coordinates $(w_{0}:\cdots :w_{n})$ 
on $
\mathbb{C}P^{n}$, we take a reduced representation $f=(f_{0}:\cdots 
:f_{n})$
which means that each $f_{i}$ is a holomorphic function on 
$\mathbb{C}^{m}$
and $f(z)=(f_{0}(z):\cdots :f_{n}(z))$ outside the analytic set $%
\{f_{0}=\cdots =f_{n}=0\}$ of codimension $\ge 2$.

Set $\Vert f\Vert := (\vert f_0\vert^2 +\cdots + \vert 
f_n\vert^2)^{1/2}.$ The characteristic function of $f$ is defined by \vskip-0.15cm
\begin{equation*}
T_{f}(r):=\int\limits_{S(r)}\log \,\Vert f\Vert \sigma
-\int\limits_{S(1)}\log \Vert f\Vert \sigma \text{ , }r>1. 
\end{equation*}
For a nonzero meromorphic function $\varphi $ on $C^{m},$ the 
characteristic
function $T_{\varphi }(r)$ of $\varphi $ is defined by considering 
$\varphi $ as a
meromorphic mapping of $\mathbb{C}^{m}$ into $\mathbb{C}P^{1}.$

Let $H=\{a_{0}w_{0}+\cdots +a_{n}w_{n}=0\}$ be a hyperplane in 
$\mathbb{C}P^{n}$
such that im$f\nsubseteq H$. Set $(f,H):=a_{0}f_{0}+\cdots 
+a_{n}f_{n}.$
We define
\begin{equation*}
N_{f}^{k)}(r,H):=N^{k)}(r,v_{(f,H)})\;\;\;\text{and}\;\;%
\;N_{l,f}^{k)}(r,H):=N_{l}^{k)}(r,v_{(f,H)}).
\end{equation*}
Sometimes we write $\overline{N}_{f}^{k)}(r,H)$ for 
$N_{1,f}^{k)}(r,H)$, $N_{l,f}(r,H)$ for $%
N_{l,f}^{+\infty )}(r,H)$ and $N_{f}(r,H)$ for  $N_{+ 
\infty,f}^{+\infty )}(r,H).$

Set $\psi_f(H) :=\displaystyle \frac{\Vert f\Vert \big( \vert a_0 
\vert^2
+\cdots + \vert a_n\vert^2\big)^{1/2}} {(f, H)}$. We define the 
proximity
function by \vskip-0.15cm
\begin{equation*}
m_{f}(r,H):=\int\limits_{S(r)}\log \,|\,\psi _{f}(H)\,|\,\sigma
-\int\limits_{S(1)}\log \,|\,\psi _{f}(H)\,|\,\sigma . 
\end{equation*}
For a nonzero meromorphic function $\varphi $, the proximity function 
is
defined by 
\begin{equation*}
m(r,\varphi ):=\int\limits_{S(r)}\log ^{+}\,|\,\ \varphi \,|\,\sigma.
\end{equation*}
We note that $m(r,\varphi)=m_{\varphi}(r,+\infty)+O(1))$ ([4], p.135).

We state the First and the Second Main Theorem of Value Distribution 
Theory.\\

\noindent
%
\textbf{First Main Theorem.}
\textit{Let} $f:\mathbb{C}^{m}\rightarrow \mathbb{C}P^{n}$ \textit{be a 
meromorphic mapping and} $H$  \textit{a hyperplane in} $\mathbb{C}%
P^{n} $ \textit{such that} im\thinspace $f\nsubseteq H$. \textit{Then 
}:\\
\begin{equation*}
N_{f}(r,H)+m_{f}(r,H)=T_{f}(r). 
\end{equation*}

\noindent\textit{For a nonzero meromorphic function }$\varphi 
$\textit{\, we have :}

\textit{\ } 
\begin{equation*}
N(r,v_{\frac{1}{\varphi }})+m(r,\varphi )=T_{\varphi }(r)+O(1). 
\end{equation*}\\

\noindent 
\noindent%
%
\textbf{Second Main Theorem.}
\textit{\textit{Let} $f:\mathbb{C}^{m}\rightarrow \mathbb{C}P^{n}$ 
\textit{be a linearly
nondegenerate meromorphic mapping and} $H_{1},...,H_{q}$ be hyperplanes 
in general position in  $\mathbb{C}P^{n}$. Then:
\begin{equation*}
(q-n-1)T_{f}(r)\le \sum_{j=1}^{q}N_{n,f}(r,H_{j})+o(T_{f}(r)) 
\end{equation*}
\textit{except for a set }$E\subset (1,+\infty )$\textit{\ of finite
Lebesgue measure.}}\\

The following so-called logarithmic derivative lemma plays an essential 
role in
Nevanlinna theory. \\

\noindent 
\noindent%
%
\textbf{Theorem 2.1.} ([5], Lemma 3.1)
\textit{Let }$\varphi $\textit{\ be a non-constant meromorphic function 
on }$%
\mathbb{C}^{m}$\textit{. Then for any }$i,\;\;1\le i\le m$\textit{,\ we 
have }

\begin{equation*}
m\Big( r,\frac{\frac{\partial }{\partial z_{i}}\varphi }{\varphi 
}\Big)%
=o(T_{\varphi }(r))\;\;\;\text{as}\;\;\;r\to \infty ,\;\;r\notin E, 
\end{equation*}
\textit{where }$E$\textit{\ }$\subset $\textit{\ }$(1,+{\infty 
})$\textit{\
of finite Lebesgue measure. }\\

Let $F, G$ and $H$ be nonzero meromorphic functions on 
$\mathbb{C}^{m}.$ For each $%
l, $ $1\le l\le m$, we define the Cartan auxiliary function by

\begin{equation*}
\Phi ^{l}(F,G,H):=F \cdot G \cdot H \cdot \left| 
\begin{array}{ccc}
1 & 1 & 1 \\ 
\frac{1}{F} & \frac{1}{G} & \frac{1}{H} \\ 
\frac{\partial }{\partial z_{l}}\left( \frac{1}{F}\right) & 
\frac{\partial }{%
\partial z_{l}}\left( \frac{1}{G}\right) & \frac{\partial }{\partial 
z_{l}}%
\left( \frac{1}{H}\right)
\end{array}
\right|. 
\end{equation*}
By [4] (Proposition 3.4) we have the following \\

\noindent 
\noindent%
%
\textbf{Theorem 2.2.}
\textit{Let $F,G,H$ be nonzero meromorphic functions on 
$\mathbb{C}^{m}$.
Assume that $\Phi ^{l}(F,G,H)\equiv 0$ and $\Phi ^{l}\Big(\displaystyle
\frac{1}{F},\displaystyle  \frac{1}{G},\displaystyle  
\frac{1}{H}\Big)\equiv
0$ for all $l$, $1\le l\le m$. Then one of the following assertions 
holds}

i) \;\;\; $F = G$ or $G = H$ or $H=F$. 

ii) $\displaystyle  \frac{F} {G}, \displaystyle  \frac{G} {H},$ 
$\displaystyle  \frac{H} {F}$ \textit{ are all constant.}

\vskip0.7cm

\section{\noindent \textbf{\ Proof of the Theorems} \vskip0.3cm}

First of all, we need the following lemmas:\\

\noindent%
%
\textbf{Lemma 1. }\textit{Let }$f_{1},...,f_{d}$\textit{\ be linearly
nondegenerate meromorphic mappings of } $\mathbb{C}^{m}$\textit{\ into 
} $%
\mathbb{C}P^{n}$\textit{\ and }$\{H_{j}\}_{j=1}^{q}$\textit{\ be
hyperplanes in } $\mathbb{C}P^{n}$\textit{\ . Then there exists a dense
subset }$%
\mathcal{C}  %
%
$\textit{\ }$\subset $\textit{\ } $\mathbb{C}^{n+1}\diagdown 
\{0\}$\textit{\
such that for any }$c=(c_{0},...,c_{n})\in 
\mathcal{C}  %
%
$\textit{, the hyperplane }$H_{c}$\textit{\ defined by }$c_{0}\omega
_{0}+...+c_{n}\omega _{n}=0$\textit{\ satisfies}
\begin{equation*}
\mathit{\ \ \ dim}(f_{i}^{-1}(H_{j})\cap f_{i}^{-1}(H_{c}))\leqslant 
m-2%
\mathit{\ for\ all\ }i\in \{1,...,d\} \textit{ and} \; j\in 
\{1,...,q\}%
\mathit{.} 
\end{equation*}\\
\textbf{Proof.} 
We refer to [5], Lemma 5.1. \hfill  $\square $\\

Let $f_{1},f_{2},f_{3}\in 
$\textit{$F_{k}\big(\{H_{j}\}_{j=1}^{q},f,1\big)$ , 
} for  $q\geq n+1.$ Set
\begin{equation*}
T(r):=T_{f_{1}}(r)+T_{f_{2}}(r)+T_{f_{3}}(r)\text{ .} 
\end{equation*}
For each $c\in 
\mathcal{C},   %
%
$ set $F_{ic}^{j}:=\displaystyle\frac{(f_{i},H_{j})}{(f_{i},H_{c})}$ 
for $i\in
\{1,2,3\}$ and $j\in \{1,\ldots ,q\}.$
\\
\noindent%
\textbf{Lemma 2.} \textit{Assume that there exist }$j_{0}\in
\{1,...,q\},c\in 
\mathcal{C}  
,l\in \{1,...,m\}$\textit{\ and a closed subset }$A$ \textit{of a 
purely }$%
(m-1)$\textit{-dimensional analytic subset of 
}$\mathbb{C}^{m}$\textit{\
satisfying}\\
\ \ \ \ \ \ \ \ 1) \qquad $\Phi _{c}^{l}:=\Phi ^{l}\big(%
F_{1c}^{j_{0}},F_{2c}^{j_{0}},F_{3c}^{j_{0}}\big)%
\not\equiv%
%
0$\textit{\  , and}\\
\ \ \ \ \ \ \ \ 2) \ \ \ \ \ \ \ 
\textit{min}$\big\{v_{(f_{1},H_{j_{0}})}^{k)},p%
\big\}=\textit{min}\big\{v_{(f_{2},H_{j_{0}})}^{k)},p\big\}=\textit{min}\big\{%
v_{(f_{3},H_{j_{0}})}^{k)},p\big\}$\textit{\ \ on 
}$\mathbb{C}^{m}\setminus
A ,$\textit{\ where }$p$\textit{\ is a positive integer. Then}
\begin{eqnarray*}
2\sum_{j=1,j\ne j_{0}}^{q}\overline{N}%
_{f_{i}}^{k)}(r,H_{j})+N_{p-1,f_{i}}^{k)}(r,H_{j_{0}}) &\leq 
&N(r,v_{\Phi
_{c}^{l}})+(p-1)\overline{N}(r,A) \\
&\leq &\frac{k+2}{k+1}T(r)+(p+2)\overline{N}(r,A)+o(T(r))
\end{eqnarray*}
\textit{\\
}$\ \ \ \ \ \ \ \ \ \ \ \ \ \ \ \ \ \ \ \ \ \ \ \ \ \ \ \ \ \ \ \ \ \ \ 
\ \
\ \ \ \ \ \ \ \ \ \ \ \ \ \ \ \ \ \ \ \ $\textit{for all }$i\in 
\{1,2,3\}$
\textit{.\\}\\
\textbf{Proof.}
Without loss of generality, we may assume that $l=1.$ For an arbitrary 
point $a\in \mathbb{C}^{m}\setminus A$ satisfying $%
v_{(f_{1},H_{j_{0}})}^{k)}(a)>0$, we have 
$v_{(f_{i},H_{j_{0}})}^{k)}(a)>0$
for all $i\in \{1,2,3\}$. We choose $a$ such that $a\notin
\bigcup\limits_{i=1}^{3}f_{i}^{-1}(H_{c})$. 
We distinguish between two cases, leading to equations (1) and (2).
\vskip0.28cm

\noindent \textbf{Case 1.} If $v_{(f_{1},H_{j_{0}})}(a)\geq p$, then $%
v_{(f_{i},H_{j_{0}})}(a)\geq p$, $i\in \{1,2,3\}$. This means that $a$ 
is a
zero point of $F_{ic}^{j_{0}}$ with multiplicity $\geq p$ for $ i\in 
\{1,2,3\}$.
We have
\begin{align*}
\Phi _{c}^{1}=& \;F_{1c}^{j_{0}}F_{3c}^{j_{0}}\frac{\partial }{\partial
z_{1}}\Big(\frac{1}{F_{3c}^{j_{0}}}\Big)-F_{1c}^{j_{0}}F_{2c}^{j_{0}}\frac{%
\partial }{\partial z_{1}}\Big(\frac{1}{F_{2c}^{j_{0}}}\Big) \\
+& \;F_{2c}^{j_{0}}F_{1c}^{j_{0}}\frac{\partial }{\partial 
z_{1}}\Big(\frac{%
1}{F_{1c}^{j_{0}}}\Big)-F_{2c}^{j_{0}}F_{3c}^{j_{0}}\frac{\partial }{%
\partial z_{1}}\Big(\frac{1}{F_{3c}^{j_{0}}}\Big) \\
+& \;F_{3c}^{j_{0}}F_{2c}^{j_{0}}\frac{\partial }{\partial 
z_{1}}\Big(\frac{%
1}{F_{2c}^{j_{0}}}\Big)-F_{3c}^{j_{0}}F_{1c}^{j_{0}}\frac{\partial }{%
\partial z_{1}}\Big(\frac{1}{F_{1c}^{j_{0}}}\Big).
\end{align*}
On the other hand 
$F_{1c}^{j_{0}}F_{3c}^{j_{0}}\displaystyle\frac{\partial 
}{\partial 
z_{1}}\Big(\frac{1}{F_{3c}^{j_{0}}}\Big)=\frac{-F_{1c}^{j_{0}}%
\frac{\partial }{\partial z_{1}}F_{3c}^{j_{0}}}{F_{3c}^{j_{0}}}$, so 
$a$ is
a zero point of 
$F_{1c}^{j_{0}}F_{3c}^{j_{0}}\displaystyle\frac{\partial }{%
\partial z_{1}}\Big(\frac{1}{F_{3c}^{j_{0}}}\Big)$ with multiplicity 
$\geq
p-1 $. By applying the same argument also to all other combinations of 
indices, we see that\\
 $a$ is a zero point of $\Phi _{c}^{1}$ with multiplicity $\geq
p-1$ .\hfill (1) \vskip0.28cm

\noindent \textbf{Case 2.} If $v_{(f_{1},H_{j_{0}})}(a)\leq p$, then $
p_{0}:=v_{(f_{1},H_{j_{0}})}(a)=v_{(f_{2},H_{j_{0}})}(a)=v_{(f_{3},H_{j_{0}})}(a)$$ 
\leq p 
$. There exists a neighborhood $U$ of $a$ such that 
$v_{(f_{1},H_{j_{0}})}%
\leq p$ on $U.$ Indeed, otherwise there exist a sequence $%
\{a_{s}\}_{s=1}^{\infty }\subset \mathbb{C}^{m}$, 
$\lim\limits_{s\rightarrow
\infty }a_{s}=a$ and $v_{(f_{1},H_{j_{0}})}(a_{s})\geq p+1$ for all 
$s$. By
the definition, we have $D^{\beta }(f_{1},H_{j_{0}})(a_{s})=0$ for all 
$%
|\beta |<p+1$. So $D^{\beta 
}(f_{1},H_{j_{0}})(a)=\lim\limits_{s\rightarrow
\infty }D^{\beta }(f_{1},H_{j_{0}})(a_{s})=0$ for all $|\beta |<p+1$. 
Thus $%
v_{(f_{1},H_{j_{0}})}(a)\geq p+1$. This is a contradiction. Hence $%
v_{(f_{1},H_{j_{0}})}\leq p$ on $U.$

We can choose $U$ such that $U\cap A=\emptyset$ , 
$v_{(f_{i},H_{j_{0}})}\leq p$
on $U$ and $(f_{i},H_{c})$ has no zero point on $U$ for all $i\in 
\{1,2,3\}$%
. Then \textit{\ 
}$v_{F_{1c}^{j_{0}}}=v_{F_{2c}^{j_{0}}}=v_{F_{3c}^{j_{0}}}%
\leq p$ on $U.$
So $U\cap \{F_{1c}^{j_{0}}=0 \}=U\cap \
\{F_{2c}^{j_{0}}=0\}=U\cap \{F_{3c}^{j_{0}}=0\}$. Choose $a$ such that 
$%
a$ is regular point of $U\cap \{F_{1c}^{j_{0}}=0\}$. By shrinking $%
U $ we may assume that there exists a holomorphic function $h$ on $U$ 
such
that $dh$ has no zero point and $F_{ic}^{j_{0}}=h^{p_{0}}u_{i}$ on $U,$
where $u_{i}(i=1,2,3)$ are nowhere vanishing holomorphic functions on 
$U$
(note that \textit{\ }$%
v_{F_{1c}^{j_{0}}}(a)=v_{F_{2c}^{j_{0}}}(a)=v_{F_{3c}^{j_{0}}}(a)=p_{0}$).
We have\vskip-0.15cm\vskip-0.15cm 
\begin{align*}
\Phi _{c}^{1}=& \;u_{1}\frac{\big(u_{3}\frac{\partial }{\partial 
z_{1}}%
u_{2}-u_{2}\frac{\partial }{\partial 
z_{1}}u_{3}\big)h^{p_{0}}}{u_{2}u_{3}}%
+u_{2}\frac{\big(u_{1}\frac{\partial }{\partial 
z_{1}}u_{3}-u_{3}\frac{%
\partial }{\partial z_{1}}u_{1}\big)h^{p_{0}}}{u_{3}u_{1}} \\
& \;\hskip4.3cm+u_{3}\frac{\big(u_{2}\frac{\partial }{\partial z_{1}}%
u_{1}-u_{1}\frac{\partial }{\partial 
z_{1}}u_{2}\big)h^{p_{0}}}{u_{1}u_{2}}.
\end{align*}
So, we have\\
$a$ is a zero point of $\Phi _{c}^{1}$ with mulitplicity $\geq p_{0}$
\hfill (2)\\

By (1), (2) and our choice of $a$, there exists an analytic set $%
M\subset \mathbb{C}^{m}$ with codimension $\geq 2$ such that $v_{\Phi
_{c}^{1}}\geq \text{min}\{v_{(f_{1},H_{j_{0}})}$ , $p-1\}$ on

\begin{equation}
\big\{z:v_{(f_{1},H_{j_{0}})}^{k)}(z)>0\big\}\setminus (M\cup A).  
\tag{3}
\end{equation}
For each $j\in \{1,\ldots ,q\}\setminus \{j_{0}\}$, let $a$ (depending 
on $j$) be an arbitrary
point in $\mathbb{C}^{m}$ such that $v_{(f_{1},H_{j})}^{k)}(a)>0$ (if 
there exist any).
 Then $
v_{(f_{i},H_{j})}^{k)}(a)>0$ for all $i\in \{1,2,3\},$ since 
$f_{1},f_{2},f_{3} \in F_{k}\big(\{H_{j}\}_{j=1}^{q},f,1\big)$. 
We can choose $a\notin
f_{i}^{-1}(H_{c})\cup f_{i}^{-1}(H_{j_{0}})$ , $i=1,2,3$. Then there 
exists
a neighborhood $U$ of $a$ such that $v_{(f_{i},H_{j})}\leq k$ on $U$ 
and $%
(f_{i},H_{j_{0}})$, $(f_{i},H_{c})$ ( $i=1,2,3$ ) have no zero point on 
$U$.
We have $B:=f_{1}^{-1}(H_{j})\cap U=f_{2}^{-1}(H_{j})\cap
U=f_{3}^{-1}(H_{j})\cap U$ and $\displaystyle\frac{1}{F_{1c}^{j_{0}}}=%
\displaystyle\frac{1}{F_{2c}^{j_{0}}}=\displaystyle\frac{1}{F_{3c}^{j_{0}}}$
on $B.$
Choose $a$ such that $a$ is a regular point of $B$. By shrinking $U,$ 
we may
assume that there exists a holomorphic function $h$ on $U$ such that 
$dh$
has no zero point and $U\cap \{h=0\}=B$. Then $\displaystyle\frac{1%
}{F_{2c}^{j_{0}}}-\displaystyle\frac{1}{F_{1c}^{j_{0}}}=h\varphi _{2}$ 
and $%
\displaystyle\frac{1}{F_{3c}^{j_{0}}}-\displaystyle\frac{1}{F_{1c}^{j_{0}}}%
=h\varphi _{3}$ on $U$ where $\varphi _{2},\varphi _{3}$ are 
holomorphic
functions on $U$. Hence, we get
\begin{eqnarray*}
\Phi _{c}^{1} &=&\;F_{1c}^{j_{0}}F_{2c}^{j_{0}}F_{3c}^{j_{0}}\left| 
\begin{array}{lll}
1 & 0 & 0 \\ 
\frac{1}{F_{1c}^{j_{0}}} & h\varphi _{2} & h\varphi _{3} \\ 
\frac{\partial }{\partial z_{1}}\Big( \frac{1}{F_{1c}^{j_{0}}}\Big) & 
\varphi _{2}\frac{\partial }{\partial z_{1}}h+h\frac{\partial 
}{\partial
z_{1}}\varphi _{2} & \,\varphi _{3}\frac{\partial }{\partial 
z_{1}}h+h\frac{%
\partial }{\partial z_{1}}\varphi _{3}
\end{array}
\right| \\
&=&F_{1c}^{j_{0}}F_{2c}^{j_{0}}F_{3c}^{j_{0}}h^{2}\left| 
\begin{array}{ll}
\varphi _{2} & \varphi _{3} \\ 
\frac{\partial }{\partial z_{1}}\varphi _{2} & \frac{\partial 
}{\partial
z_{1}}\varphi _{3}
\end{array}
\right|.
\end{eqnarray*}
Hence,  $a$ is a zero point of $\Phi _{c}^{1}$ with multiplicity $\ge 
2$. Thus,
for each $j\in \{1,\ldots ,q\}\setminus \{j_{0}\}$, there exists an 
analytic
set $N\subset \mathbb{C}^{m}$ with codimension $\ge 2$ such that 
$v_{\Phi
_{c}^{1}}\ge 2$ on
\begin{equation}
\big\{z:v_{(f_{1},H_{j})}^{k)}(z)>0\big\}\setminus N.  \tag{4}
\end{equation}
By (3) and (4), we have
\begin{equation*}
\ 2\sum_{j=1,j\ne j_{0}}^{q}\overline{N}_{f_{1}}^{k)}(r,H_{j})+N
_{p-1,f_{1}}^{k)}(r,H_{j_{0}})\leq N(r,v_{\Phi 
_{c}^{1}})+(p-1)\overline{N}%
(r,A). 
\end{equation*}
Similarly, we have
\begin{equation}
2\sum_{j=1,j\ne j_{0}}^{q}\overline{N\text{ }}%
_{f_{i}}^{k)}(r,H_{j})+N_{p-1,f_{i}}^{k)}(r,H_{j_{0}})\leq N(r,v_{\Phi
_{c}^{1}})+(p-1)\overline{N}(r,A),\text{ \ }i=1,2,3.\text{\ \ \ \ \ \ \ 
} 
\tag{5}
\end{equation}
Let $a$ be an arbitrary zero point of some $F_{ic}^{j_{0}}$ $,$ 
$a\notin
f_{i}^{-1}(H_{c})$, say $i=1.$
We have\vskip-0.15cm\vskip-0.15cm 
\begin{align}
\Phi _{c}^{1}=& \;\big( F_{2c}^{j_{0}}-F_{3c}^{j_{0}}\big) 
F_{1c}^{j_{0}}%
\frac{\partial }{\partial z_{1}}\Big( \frac{1}{F_{1c}^{j_{0}}}\Big) 
+\big( %
F_{3c}^{j_{0}}-F_{1c}^{j_{0}}\big) F_{2c}^{j_{0}}\frac{\partial 
}{\partial
z_{1}}\Big( \frac{1}{F_{2c}^{j_{0}}}\Big)  \notag \\
+& \;\big( F_{1c}^{j_{0}}-F_{2c}^{j_{0}}\big) 
F_{3c}^{j_{0}}\frac{\partial }{%
\partial z_{1}}\Big( \frac{1}{F_{3c}^{j_{0}}}\Big).  \tag{6}
\end{align}
So we have \vskip-0.15cm\vskip-0.15cm 
\begin{equation*}
v_{\frac{1}{\Phi _{c}^{1}}}(a)\le 1+\max \{v_{\frac{1}{F_{ic}^{j_{0}}}%
}(a),i=2,3\}\;\leq 1+\text{ 
}v_{\frac{1}{F_{2c}^{j_{0}}}}(a)\;+v_{\frac{1}{%
F_{3c}^{j_{0}}}}(a)\text{ .}\; 
\end{equation*}
Furthermore, if $0<v_{F_{1c}^{j_{0}}}(a)\le k$ \big(and, hence, $ 
v_{(f_{1},H_{j_{0}})}^{k)}(a)>0\big)$ and $a\notin A$, then by (3) we may
assume that $v_{\frac{1}{\Phi _{c}^{1}}}(a)=0$ (outside an analytic set
of codimension $\geq 2$ ). \hfill(7) \\
Let $a$ be an arbitrary pole of all $F_{ic}^{j_{0}}$ , $i=1,2,3.$ By 
(6) we
have 
\begin{equation}
v_{\frac{1}{\Phi _{c}^{1}}}(a)\le \max \{v_{\frac{1}{F_{ic}^{j_{0}}}%
}(a),i=1,2,3\}+1<\sum_{i=1}^{3}v_{\frac{1}{F_{ic}^{j_{0}}}}(a)  \tag{8}
\end{equation}
It follows from (6) that a pole of $\Phi _{c}^{1}$ is a zero or a pole 
of some $%
F_{ic}^{j_{0}}.$ Thus, by (6), (7) and (8), we have
\begin{align}
N\Big(r,v_{\frac{1}{\Phi _{c}^{1}}}\Big) \le & \;\sum_{i=1}^{3}N\Big( 
r,v_{%
\frac{1}{F_{ic}^{j_{0}}}}\Big) +\sum_{i=1}^{3}\Big(\overline{N}\big(%
r,v_{F_{ic}^{j_{0}}}\big) 
-\overline{N}^{k)}\big(r,v_{F_{ic}^{j_{0}}}\big) %
\Big)+3\overline{N}(r,A)  \notag \\
\le & \;\sum_{i=1}^{3}N\Big( 
r,v_{\frac{1}{F_{ic}^{j_{0}}}}\Big)+\frac{1}{k+1%
}\sum_{i=1}^{3}N\big( r,v_{F_{ic}^{j_{0}}}\big) +3\overline{N}(r,A) 
\notag \\
& \le \sum_{i=1}^{3}N\big( 
r,v_{\frac{1}{F_{ic}^{j_{0}}}}\big)+\frac{1}{k+1}%
\sum_{i=1}^{3}T_{F_{ic}^{j_{0}}}(r)+3\overline{N}(r,A)  \notag \\
& \le \sum_{i=1}^{3}N\big( 
r,v_{\frac{1}{F_{ic}^{j_{0}}}}\big)+\frac{1}{k+1}%
T(r)+3\overline{N}(r,A)+O(1).  \tag{9}
\end{align}
We have \vskip-0.15cm\vskip-0.15cm
\begin{align*}
\Phi_c^1 =&\; F_{1c}^{j_0} \Big[F_{3c}^{j_0} \frac{\partial} {\partial 
z_1} %
\Big( \frac{1} {F_{3c}^{j_0}}\Big)- F_{2c}^{j_0} \frac{\partial} 
{\partial
z_1} \Big( \frac{1} {F_{2c}^{j_0}}\Big)\Big] \\
+ &\; F_{2c}^{j_0} \Big[F_{1c}^{j_0} \frac{\partial} {\partial z_1} 
\Big( 
\frac{1} {F_{1c}^{j_0}}\Big)- F_{3c}^{j_0} \frac{\partial} {\partial 
z_1} %
\Big( \frac{1} {F_{3c}^{j_0}}\Big)\Big] \\
+&\; F_{3c}^{j_0} \Big[F_{2c}^{j_0} \frac{\partial} {\partial z_1} 
\Big( 
\frac{1} {F_{2c}^{j_0}}\Big)- F_{1c}^{j_0} \frac{\partial} {\partial 
z_1} %
\Big( \frac{1} {F_{1c}^{j_0}}\Big)\Big]
\end{align*}
so $m (r, \Phi_c^1) \le \sum\limits_{i=1}^3 m (r, F_{ic}^{j_0}) + 2
\sum\limits_{i=1}^3 m\Big( r, F_{ic}^{j_0} \frac{\partial} {\partial 
z_1} %
\Big(\frac{1} {F_{ic}^{j_0}}\Big) \Big) + 0 (1).$
By Theorem 2.1, we have
\begin{equation*}
m\Big( r,F_{ic}^{j_{0}}\frac{\partial }{\partial z_{1}}\Big(\frac{1}{%
F_{ic}^{j_{0}}}\Big) \Big) =o\big( T_{F_{ic}^{j_{0}}}(r)\big)
\end{equation*}
Thus, we get \vskip-0.15cm\vskip-0.15cm
\begin{equation}
m(r,\Phi _{c}^{1})\le \sum_{i=1}^{3}m(r,F_{ic}^{j_{0}})+o(T(r)),  
\tag{10}
\end{equation}
(note that $T_{F_{ic}^{j_{0}}}(r)\leq T_{f_{i}}(r)\,+O(1)$).\\
By (9) , (10) and by the First Main Theorem, we have
\begin{eqnarray*}
N(r,v_{\Phi _{c}^{1}}) &\leq & T_{\Phi 
_{c}^{1}}(r)+O(1)=N\big(r,v_{\frac{1%
}{\Phi _{c}^{1}}}\big)+m(r,\Phi _{c}^{1})+O(1) \\
&\leq &\;\sum_{i=1}^{3}\Big(N\big(r,v_{\frac{1}{F_{ic}^{j_{0}}}}\big)%
+m(r,F_{ic}^{j_{0}})\Big)+\frac{1}{k+1}T(r)+3\overset{\_\_}{N}(r,A)+o(T(r))
\\
&\leq 
&\;\sum_{i=1}^{3}T_{F_{ic}^{j_{0}}}(r)+\frac{1}{k+1}T(r)+3\overset{\_\_%
}{N}(r,A)+o(T(r)) \\
&\leq 
&\;\sum_{i=1}^{3}T_{f_{i}}(r)+\frac{1}{k+1}T(r)+3\overset{\_\_}{N}%
(r,A)+o(T(r)) \\
&=&\frac{k+2}{k+1}T(r)+3\overset{\_\_}{N}(r,A)+o(T(r)). \hspace{4.3cm} 
\text{ (11) }
\end{eqnarray*}
By (5) and (11) we get Lemma 2. \hfill $\square $\\

The following lemma is a variant of the Second Main Theorem  
without taking account of multiplicities of order $>k$ in the counting 
functions.\\

\noindent%
%
\textbf{Lemma 3.} \textit{Let }$f$\textit{\ be a linearly nondegenerate
meromorphic mapping of }$\mathbb{C}^{m}$\textit{\ into } 
$\mathbb{C}P^{n}$%
\textit{\ and }$\{H_{j}\}_{j=1}^{q}(q\geq n+2)$\textit{\ be hyperplanes 
in } 
$\mathbb{C}P^{n}$\textit{\ in general position. Take a positive integer 
}$k$%
\textit{ with }$\frac{qn}{q-n-1}\leq k\leq +\infty $\textit{\ . Then}
\begin{eqnarray*}
\mathit{T}_{f}\mathit{(r)} &\leq &\frac{k}{(q-n-1)(k+1)-qn}%
\sum_{j=1}^{q}N_{n,f}^{k)}(r,H_{j})\mathit{+o(T}_{f}\mathit{(r))} \\
&\leq 
&\frac{nk}{(q-n-1)(k+1)-qn}\sum_{j=1}^{q}\overline{N}_{f}^{k)}\mathit{%
(r,H}_{j}\mathit{)+o(T}_{f}\mathit{(r))}
\end{eqnarray*}
\textit{for all }$r>1$\textit{\ except a set }$E$\textit{\ of finite
Lebesgue measure.}\\

\noindent\textbf{Proof.}
By the First and the Second Main Theorems, we have 
\begin{equation*}
(q-n-1)T_{f}(r)\leq \sum_{j=1}^{q}N_{n,f}(r,H_{j})+o\big(T_{f}(r)\big)
\end{equation*}
\begin{align*}
& \leq \frac{k}{k+1}\sum_{j=1}^{q}N_{n,f}^{k)}(r,H_{j})+\frac{n}{k+1}%
\sum_{j=1}^{q}N_{f}(r,H_{j})+o\big(T_{f}(r)\big) \\
& \leq \frac{k}{k+1}\sum_{j=1}^{q}N_{n,f}^{k)}(r,H_{j})+\frac{qn}{k+1}%
T_{f}(r)+o\big(T_{f}(r)\big),\;\;r\notin E,
\end{align*}
which impies that
\begin{equation*}
\Big(q-n-1-\frac{qn}{k+1}\Big)T_{f}(r)\leq \frac{k}{k+1}%
\sum_{j=1}^{q}N_{n,f}^{k)}(r,H_{j})+o\big(T_{f}(r)\big)
\text{ .} 
\end{equation*}
Thus, we have
\begin{eqnarray*}
\mathit{T}_{f}\mathit{(r)} &\leq &\frac{k}{(q-n-1)(k+1)-qn}%
\sum_{j=1}^{q}N_{n,f}^{k)}(r,H_{j})\mathit{+o(T}_{f}\mathit{(r))} \\
&\leq 
&\frac{nk}{(q-n-1)(k+1)-qn}\sum_{j=1}^{q}\overline{N}_{f}^{k)}\mathit{%
(r,H}_{j}\mathit{)+o(T}_{f}\mathit{(r))}\text{ } \;\;\;\;\;\; \square
\end{eqnarray*} \\

\noindent \textbf{Proof of Theorem 1.}
Assume that \ there exist three distinct mappings $f_{1},f_{2},f_{3}\in 
$\ $%
F_{k}(\{H_{j}\}_{j=1}^{q},f,p).$
Denote by $Q$ the set which contains all indices $j\in \{1,...,q\}$ 
satisfying 
$\Phi ^{l}\big(F_{1c}^{j},F_{2c}^{j},F_{3c}^{j}\big)%
\not\equiv%
0$\textit{\ }for some $c\in 
\mathcal{C}  
$\textit{\ } and some $l\in \{1,...,m\}.$
We now prove that 
\begin{equation}
\#(\{1,...,q\}\backslash Q)\geq
3n-1. \tag{12}
\end{equation}
For the proof of (12) we distinguish three cases:\\

\noindent \textbf{Case 1.} $1\leq n\leq 3,q=3n+1,p=2,k\geq 23n$ $.$\\
Suppose that (12) does not hold, then $\#Q\geq 3.$
For each $j_{0}\in Q$, by Lemma 2 (with $A=$\ $\emptyset$, $p=2$) we 
have
\begin{equation}
\ \ \ 2\sum_{j=1,j\ne 
i_{0}}^{q}\overline{N}_{f_{i}}^{k)}(r,H_{j})+\overline{%
N}_{f_{i}}^{k)}(r,H_{j_{0}})\leq \frac{k+2}{k+1}T(r)+o(T(r))\ ,\ 
i=1,2,3.\ \
\   \tag{13}
\end{equation}
By (13) and Lemma 3 we have
\begin{align*}
& \Big(q-n-1-\frac{qn}{k+1}\Big)T_{f_{i}}(r)\leq 
\frac{nk}{k+1}\sum_{j=1}^{q}%
\overline{N}_{f_{i}}^{k)}(r,H_{j})+o\big(T_{f_{i}}(r)\big) \\
& \leq \;\frac{nk(k+2)}{2(k+1)^{2}}T(r)+\frac{nk}{2(k+1)}%
\overline{N}_{f_{i}}^{k)}(r,H_{j_{0}})+o\big(T(r)\big),i=1,2,3.
\end{align*}
Thus, we obtain\vskip-0.15cm\vskip-0.15cm 
\begin{align*}
 \Big(q-n-1-\frac{qn}{k+1}\Big)T(r)\leq 
\frac{3nk(k+2)}{2(k+1)^{2}}T(r)+%
\frac{nk}{2(k+1)}\sum_{i=1}^{3}\overline{N}_{f_{i}}^{k)}(r,H_{j_{0}})+o\big(%
T(r)\big),
\end{align*}
which implies
\begin{align*}
 \big[2(q-n-1)(k+1)^{2}-2qn(k+1)-3nk(k+2)\big]T(r)  \\
 \leq nk(k+1)\sum_{i=1}^{3}\overline{N}
_{f_{i}}^{k)}(r,H_{j_{0}})+o(T(r))=3nk(k+1)\overline{N}
_{f_{i}}^{k)}(r,H_{j_{0}})+o(T(r))\text{ .}
\end{align*}
Hence, we have
\begin{align*}
 \liminf\limits_{r\rightarrow \infty \;\,r\not{\in}E}\frac{
\overline{N}_{f_{i}}^{k)}(r,H_{j_{0}})}{T(r)}&\geq \frac{
2(q-n-1)(k+1)^{2}-2qn(k+1)-3nk(k+2)}{3nk(k+1)}\\ 
&=\frac{k^{2}-6nk-6n+2}{3k(k+1)}\text{  , }i\in\{1,2,3\}.\tag{14}
\end{align*}
Set \vskip-0.15cm\vskip-0.15cm
\begin{equation*}
A_{i}:=\big\{ r>1:T_{f_{i}}(r)=\textit{min}\big\{ %
T_{f_{1}}(r),T_{f_{2}}(r),T_{f_{3}}(r)\big\}\big\},\;\;\;i\in 
\{1,2,3\}. 
\end{equation*}
Then $A_{1}\cup A_{2}\cup A_{3}=(1,+{\infty }).$ Without loss of 
generality, we may 
assume that the Lebesgue measure of $A_{1}$ is infinite.
By (14) we have
\begin{equation*}
\liminf\limits_{r\rightarrow \infty \;\,r\in A_{1}\setminus E}\frac{%
\overline{N}_{f_{1}}^{k)}(r,H_{j_{0}})}{T_{f_{1}}(r)}\geq \frac{%
k^{2}-6nk-6n+2}{k(k+1)}\text{ },j_{0}\in Q. 
\end{equation*}

 Take three distinct indices $j_{1},j_{2},j_{3}\in Q$ $($note that $%
\#Q\geq 3)$. Then we have
\begin{align*}
& \liminf\limits_{r\to \infty \;\,r\in A_{1}\setminus 
E}\frac{\overline{N}%
_{f_{1}}^{k)}(r,H_{j_{1}})+\overline{N}_{f_{1}}^{k)}(r,H_{j_{2}})+\overline{N%
}_{f_{1}}^{k)}(r,H_{j_{3}})}{T_{f_{1}}(r)}\ge 
\frac{3(k^{2}-6nk-6n+2)}{k(k+1)%
},  
\end{align*}
which implies that
\begin{align*}
\liminf\limits_{r\to \infty \;\,r\in A_{1}\setminus E}\frac{%
\sum\limits_{j=1}^{q}\overline{N}_{f_{1}}^{k)}(r,H_{j})}{T_{f_{1}}(r)}\ge 
\frac{3(k^{2}-6nk-6n+2)}{k(k+1)}\text{ .}  \tag{15}
\end{align*}

Since $f_{1}%
\not\equiv%
f_{2}$ there exists $c\in 
\mathcal{C}  
$\textit{\ }such that\textit{\ } $\frac{(f_{1},H_{1})}{(f_{1},H_{c})}$ 
$%
\not\equiv%
\frac{(f_{2},H_{1})}{(f_{2},H_{c})}$. Indeed, otherwise by Lemma 1 we 
have
that $\frac{(f_{1},H_{1})}{(f_{1},H)}$ $\equiv 
\frac{(f_{2},H_{1})}{(f_{2},H)%
}$ for all hyperplanes $H$ in $\mathbb{C}P^{n}$ . In particular 
$\frac{%
(f_{1},H_{1})}{(f_{1},H_{j})}$ $\equiv 
\frac{(f_{2},H_{1})}{(f_{2},H_{j})}$
for all \ $j=2,...,n+1.$
We choose homogeneous coordinates $(\omega _{0}:\cdots :\omega _{n})$ 
on $\mathbb{C}%
P^{n}$ with $H_{j}=\{\omega _{j}=0\}$ $(1\le j\le n+1)$
 and take reduced representations: $f_{1}=(f_{1_{1}}:\cdots
:f_{1_{n+1}})$, $f_{2}=(f_{2_{1}}:\cdots :f_{2_{n+1}}).$ Then
\begin{equation*}
\begin{cases}\displaystyle\frac{f_{1_{j}}}{f_{1_{1}}}=\displaystyle\frac{f_{2_{j}}}{%
f_{2_{1}}}\\
\text{ (}j=2,\ldots ,n+1)\end{cases}\;\;\Rightarrow 
\;\;\frac{f_{1_{1}}}{f_{2_{1}}}%
=\cdots =\frac{f_{1_{n+1}}}{f_{2_{n+1}}}\;\Rightarrow f_{1}\equiv 
f_{2}. 
\end{equation*}
This is a contradiction.

Since \textit{dim} $(f_{i}^{-1}(H_{1})\cap f_{i}^{-1}(H_{c}))\leq m-2$ 
we have
$$
\begin{array}{lll}
\ T_{\frac{(f_{i},H_{1})}{(f_{i},H_{c})}}(r) & = & 
\int\limits_{S(r)}\log
\,(\left| (f_{i},H_{1})\right| ^{2}+\left| (f_{i},H_{c})\right| 
^{2})^{\frac{%
1}{2}}\sigma \ \ +O(1) \\ 
& \leq & \int\limits_{S(r)}\log \,\Vert f_{i}\Vert \sigma
+O(1)=T_{f_{i}}(r)+O(1)\ ,\; i=1,2,3.
\end{array}
$$
Since\ $f_{1}=f_{2}$ on $\bigcup\limits_{j=1}^{q}\big\{%
z:v_{(f_{1},H_{j})}^{k)}(z)>0\big\}$ and\\
 \textit{dim}$\Big\{%
z:v_{(f_{1},H_{i})}^{k)}(z)>0\ \ and\ \ 
v_{(f_{1},H_{j})}^{k)}(z)>0\Big\}%
\leq m-2$ for all $i\neq j,$ ? ?we have\\
 \begin{align*}
& \sum_{j=1}^{q}\overline{N}_{f_{1}}^{k)}(r,H_{j})\leq 
N\big(r,v_{\frac{%
(f_{1},H_{1})}{(f_{1},H_{c})}-\frac{(f_{2},H_{1})}{(f_{2},H_{c})}}\big)\leq
T_{\frac{(f_{1},H_{1})}{(f_{1},H_{c})}-\frac{(f_{2},H_{1})}{(f_{2},H_{c})}%
}(r)\text{ }+0(1) \\
& \leq 
T_{\frac{(f_{1},H_{1})}{(f_{1},H_{c})}}(r)+T_{\frac{(f_{2},H_{1})}{%
(f_{2},H_{c})}}(r)+0(1)\leq T_{f_{1}}(r)+T_{f_{2}}(r)+0(1),\; 
\end{align*}
which implies
\begin{align*}
 \liminf\limits_{r\rightarrow \infty }\frac{%
T_{f_{1}}(r)+T_{f_{2}}(r)}{\sum\limits_{j=1}^{q}\overline{N}%
_{f_{1}}^{k)}(r,H_{j})}\geq 1.
\end{align*}
On the other hand, by Lemma 3, we have
\begin{align*}
& \Big(q-n-1-\frac{qn}{k+1}\Big) T_{f_{i}}(r)\le 
\frac{nk}{k+1}\sum_{j=1}^{q}%
\overline{N}_{f_{i}}^{k)}(r,H_{j})+o\big(T_{f_{i}}(r)\big) \\
& 
\hskip4cm=\frac{nk}{k+1}\sum_{j=1}^{q}\overline{N}_{f_{1}}^{k)}(r,H_{j})+o%
\big( T_{f_{i}}(r)\big)\text{ ,}
\end{align*}
which implies
\begin{equation*}
 \limsup\limits_{r\to \infty \;\,r\not{\in}E}\frac{T_{f_{i}}(r)}{%
\sum\limits_{j=1}^{q}\overline{N}_{f_{1}}^{k)}(r,H_{j})}\le \frac{nk}{%
(q-n-1)(k+1)-qn}\text{ , }i=1,2,3\text{ .} 
\end{equation*}
Hence, we obtain \vskip-0.15cm\vskip-0.15cm 
$$
\limsup\limits_{r\to \infty \;\,r\in A_{1}\setminus 
E}\frac{T_{f_{1}}(r)}{%
\sum\limits_{j=1}^{q}\overline{N}_{f_{1}}^{k)}(r,H_{j})}
=
\limsup\limits_{r\to \infty \;\,r\in A_{1}\setminus 
E}\big(\frac{T_{f_{1}}(r)+
T_{f_{2}}(r)}{
\sum\limits_{j=1}^{q}\overline{N}_{f_{1}}^{k)}(r,H_{j})}-
\frac{T_{f_{2}}(r)}{
\sum\limits_{j=1}^{q}\overline{N}_{f_{1}}^{k)}(r,H_{j})} \big)$$
$$\geq
\liminf\limits_{r\to \infty \;\,r\in A_{1}\setminus 
E}\frac{T_{f_{1}}(r)+
T_{f_{2}}(r)}{
\sum\limits_{j=1}^{q}\overline{N}_{f_{1}}^{k)}(r,H_{j})}\;
-
\limsup\limits_{r\to \infty \;\,r\in A_{1}\setminus 
E}\frac{T_{f_{2}}(r)}{
\sum\limits_{j=1}^{q}\overline{N}_{f_{1}}^{k)}(r,H_{j})}
\ge
 1-\frac{nk}{%
(q-n-1)(k+1)-qn} 
$$
 Consequently, we get
\begin{eqnarray*}
\liminf\limits_{r\to \infty \;\,r\in A_{1}\setminus E}\frac{%
\sum\limits_{j=1}^{q}\overline{N}_{f_{1}}^{k)}(r,H_{j})}{T_{f_{1}}(r)} 
&\le &%
\frac{(q-n-1)(k+1)-qn}{(q-n-1)(k+1)-qn-nk} \\
&=&\frac{2k+1-3n}{k+1-3n} \hspace{3cm}(16) ¾
\end{eqnarray*}

By (15) and (16) we have \vskip-0.15cm\vskip-0.15cm 
\begin{equation*}
\frac{3(k^{2}-6nk-6n+2)}{k(k+1)}\leq \frac{2k+1-3n}{k+1-3n}. 
\end{equation*}
This contradicts $k\geq 23n$. Thus, we get (12) in this case.\\

\noindent \textbf{Case 2. \ \ }$4\leq n\leq 6$, $q=3n$, $p=2$, 
$k\geq \frac{(6n-1)n}{n-3}$ \ .\\
Suppose that (12) does not hold, then there exists $j_{0}\in Q.$
By Lemma 2 (with $A=\emptyset,p=2$) we have
\begin{equation*}
\ \ \ \ \ \ 2\sum_{j=1,j\ne 
j_{0}}^{3n}\overline{N}_{f_{i}}^{k)}(r,H_{j})+%
\overline{N}_{f_{i}}^{k)}(r,H_{j_{0}})\leq \frac{k+2}{k+1}T(r)+o(T(r))\ 
,\
i=1,2,3\ . 
\end{equation*}
On the other hand, by Lemma 3 we have
\begin{eqnarray*}
\text{ }\sum_{j=1,j\ne j_{0}}^{3n}\overline{N}%
_{f_{i}}^{k)}(r,H_{j})+o(T_{f_{i}}(r)) &\geq 
&\frac{(2n-2)(k+1)-(3n-1)n}{nk}%
T_{f_{i}}(r), \text{ and} \\
\text{ 
}\sum_{j=1}^{3n}\overline{N}_{f_{i}}^{k)}(r,H_{j})+o(T_{f_{i}}(r))
&\geq &\frac{(2n-1)(k+1)-3n^{2}}{nk}T_{f_{i}}(r),
\end{eqnarray*}
which implies that
\begin{equation*}
 \text{ }2\sum_{j=1,j\ne j_{0}}^{3n}\overline{N}%
_{f_{i}}^{k)}(r,H_{j})+\overline{N}%
_{f_{i}}^{k)}(r,H_{j_{0}})+o(T_{f_{i}}(r))\geq 
\frac{(4n-3)(k+1)-(6n-1)n}{nk}%
T_{f_{i}}(r) 
\end{equation*}
Hence, we have 
\begin{equation*}
\ \frac{(4n-3)(k+1)-(6n-1)n}{nk}T_{f_{i}}(r)\leq \frac{k+2}{k+1}%
T(r)+o(T(r))\ ,\ i=1,2,3. 
\end{equation*}
Consequently, we get
\begin{align*}
 \ \frac{(4n-3)(k+1)-(6n-1)n}{nk}T(r) \leq \frac{3(k+2)}{k+1}%
T(r)+o(T(r)), 
\end{align*}
which implies that
\begin{align*}
  \ \left( (4n-3)(k+1)-(6n-1)n\right) T(r) & \leq & \frac{%
3nk(k+2)}{k+1}T(r)+o(T(r)) \\ 
&   \leq & 3n(k+1)T(r)+o(T(r). 
\end{align*}
Hence, we obtain
$k+1\leq \frac{(6n-1)n}{n-3}. $   \ This is a contradiction. Thus, we 
get (12) in this case$.$\\

\noindent \textbf{Case 3. } $n\geq 7$, $q=3n-1$, $p=1$, $k\geq 
\frac{(6n-4)n}
{n-6}$ $.$\\
Suppose that (12) does not hold, then there exists $j_{0}\in Q.$
By Lemma 2 (with $A=\emptyset , p=1$) we have
\begin{equation*}
\ \ \ 2\text{ }\sum_{j=1,j\ne j_{0}}^{3n-1}\overline{N}%
_{f_{i}}^{k)}(r,H_{j})\leq \frac{k+2}{k+1}T(r)+o(T(r))\ ,\ i=1,2,3\  
\end{equation*}
 (note that $%
N_{0,f_{i}}^{k)}(r,H_{j_{0}})=0$).
On the other hand, by Lemma 3, we have
\begin{equation*}
\ \ 2\sum_{j=1,j\ne j_{0}}^{3n-1}\overline{N}%
_{f_{i}}^{k)}(r,H_{j})+o(T_{f_{i}}(r))\geq 
2\frac{(2n-3)(k+1)-(3n-2)n}{nk}%
T_{f_{i}}(r) 
\end{equation*}
 Hence, we get
\begin{align*}
 \frac{2[(2n-3)(k+1)-(3n-2)n)]}{nk}T_{f_{i}}(r)  \leq  \frac{k+2}{k+1}%
T(r)+o(T(r))\ ,
 \end{align*}
which implies
\begin{align*}
  ((4n-6)(k+1)-(6n-4)n)T_{f_{i}}(r)  \leq \frac{nk(k+2)}{k+1}%
T(r)+o(T(r)),\ i=1,2,3.
\end{align*}
Hence, we have
\begin{align*}
  ((4n-6)(k+1)-(6n-4)n)T(r)  &\leq&  \frac{3nk(k+2)}{k+1}
T(r)+o(T(r)) \\ 
   &\leq & 3n(k+1)T(r)+o(T(r)). 
\end{align*}
Thus, we obtain 
$$(4n-6)(k+1)-(6n-4)n  \leq  3n(k+1)$$
implying
$$´k+1\leq \frac{(6n-4)n}{n-6},$$ 
which is a contradiction. Thus, we get (12) in this case.\\

So, for any case we have $\#(\{1,\ldots ,q\}\setminus Q)\geq 3n-1$. 
Without
loss of generality, we may assume that $1,\ldots ,3n-1\notin Q.$
We have

$\Phi ^{l}\big(F_{1c}^{j},F_{2c}^{j},F_{3c}^{j}\big)\equiv 0$ for all 
$c\in 
\mathcal{C}  %
%
$, $l\in \{1,...,m\},$ $j\in \{1,\ldots ,3n-1\}$.

\noindent On the other hand, $%
\mathcal{C}  %
%
$ is dense in $\mathbb{C}^{n+1}$ . Hence, $\Phi ^{l}\big(%
F_{1c}^{j},F_{2c}^{j},F_{3c}^{j}\big )\equiv 0$ for all $c\in 
\mathbb{C}%
^{n+1}\setminus \{0\}$, $l\in \{1,...,m\},$ $j\in \{1,\ldots ,3n-1\}$. 
In
particular (for $H_{c}=H_{i})$ we have
\begin{equation*}
\Phi ^{l}\left( 
\frac{(f_{1},H_{j})}{(f_{1},H_{i})},\;\;\frac{(f_{2},H_{j})}{%
(f_{2},H_{i})},\;\;\frac{(f_{3},H_{j})}{(f_{3},H_{i})}\right) \equiv 0 
\end{equation*}
for all $1\le i\ne j\le 3n-1,\;\;\;l\in \{1,\ldots ,m\}.$\hfill (17)\\

In the following we distinguish between the cases $n=1$ and $n \geq 
2$.\\

\noindent \textbf{Case 1.} If $n=1$, then $a_{j}:=H_{j}(j=1,2,3,4)$ are 
distinct points in 
$\mathbb{C}P^{1}$.
We have that
\begin{equation*}
g_{1}:=\frac{(f_{1},a_{1})}{(f_{1},a_{2})},\;\;g_{2}:=\frac{(f_{2},a_{1})}{%
(f_{2},a_{2})},\;\;g_{3}:=\frac{(f_{3},a_{1})}{(f_{3},a_{2})} 
\end{equation*}
are distinct nonconstant meromorphic functions. 
By (17) and by Theorem 2.2, there exist constants $\alpha ,\beta $ such 
that
\begin{align*}
g_{2}=\alpha g_{1\text{ \ \ }}
\text{ , \ }g_{3}=\beta g_{1}\text{ },(\alpha ,\beta \notin \{1,\infty
,0\},\alpha \neq \beta )  \tag{18}
\end{align*}
We have $v_{(f_1, a_3)} \ge k+1$ on $\{z: (f_1, a_3) (z) =0\}$:
Indeed, otherwise there exists $z_{0}$ such that 
$0<v_{(f_{1},a_{3})}(z_{0})%
\le k$. Then $v_{(f_{i},a_{3})}^{k)}(z_{0})>0$, for all $i\in 
\{1,2,3\}$. We
have $(f_{1},a_{3})(z_{0})=(f_{2},a_{3})(z_{0})=0$ $\Rightarrow
\;f_{1}(z_{0})=f_{2}(z_{0})=a_{3}^{*},$ where we denote $%
a_{j}^{*}:=(a_{j_{1}}:-a_{j_{0}})$ for every point $%
a_{j}=(a_{j_{0}}:a_{j_{1}})\in \mathbb{C}P^{1}$. So $%
g_{1}(z_{0})=g_{2}(z_{0})=\displaystyle\frac{(a_{3}^{*},a_{1})}{%
(a_{3}^{*},a_{2})}\ne 0$, $\infty $ (note that $a_{3}\ne a_{1}$, 
$a_{3}\ne
a_{2})$. So, by (18) we have $\alpha =1$. This is a contradiction. 
Thus $v_{(f_{1},a_{3})}\ge k+1$ on $\{z:(f_{1},a_{3})(z)=0\}$. 
Similarly, $%
v_{(f_{i},a_{j})}\ge k+1$ on $\{z:(f_{i},a_{j})(z)=0\}$ for $i\in 
\{1,2,3\}$%
, $j\in \{3,4\}$.

Set $b_1 =\alpha\displaystyle \frac{(a_3^*, a_1)} {(a_3^*, 
a_2)},\;\;b_2 =%
\displaystyle\frac{\alpha}{\beta}\displaystyle \frac{(a_3^*, a_1)} 
{(a_3^*,
a_2)},\;\; b_3 =\displaystyle \frac{(a_3^*, a_1)} {(a_3^*, a_2)}.$
Then we have 
\begin{align*}
v_{g_{2}-b_{3}}&=v_{\frac{_{(f_{2},a_{3})(a_{1}^{*},a_{2})}}{%
(f_{2},a_{2})(a_{3}^{*},a_{2})}}\;\geq k+1\;\;\text{on}\;\;%
\{z:(g_{2}-b_{3})(z)=0\}, \\
v_{g_{2}-b_{1}}= \;v_{g_{1}-\frac{1}{\alpha }b_{1}}&=v_{\frac{%
_{(f_{1},a_{3})(a_{1}^{*},a_{2})}}{(f_{1},a_{2})(a_{3}^{*},a_{2})}}\geq
k+1\;\;\text{on}\;\;\{z:(g_{2}-b_{1})(z)=0\},\text{and} \\
v_{g_{2}-b_{2}}= \;v_{g_{3}-\frac{\beta }{\alpha }b_{2}}&=v_{\frac{%
_{(f_{3},a_{3})(a_{1}^{*},a_{2})}}{(f_{3},a_{2})(a_{3}^{*},a_{2})}}\geq
k+1\;\;\text{on}\;\;\{z:(g_{2}-b_{2})(z)=0\}.\;\;
\end{align*}
Since the points $b_{1},b_{2},b_{3}$ are distinct,
by the First and the Second Main Theorem, we have
\begin{align*}
T_{g_{2}}(r)& \leq 
\sum_{j=1}^{3}\overline{N}\big(r,v_{g_{2}-b_{j}}\big)+o%
\big(T_{g_{2}}(r)\big) \\
& \leq \frac{1}{k+1}\sum_{j=1}^{3}N\big(r,v_{g_{2}-b_{j}}\big)+o\big(%
T_{g_{2}}(r)\big) \\
& \leq \frac{3}{k+1}T_{g_{2}}(r)+o\big(T_{g_{2}}(r)\big)\text{ .}
\end{align*}
This contradicts $k\geq 23.$\\

\noindent \textbf{Case 2.} If $n\ge 2$ ,
for each $1\le i\ne j\le 3n-1$, by (17) and Theorem 2.2., there exists 
a
constant $\alpha _{ij}$ such that
\begin{equation*}
\frac{(f_2, H_j)} {(f_2, H_i)} =\alpha_{ij} \frac{(f_1, H_j)} {(f_1, 
H_i)}
\text{ or }\frac{(f_3, H_j)} {(f_3, H_i)}=\alpha_{ij} \frac{(f_1, H_j)} 
{(f_1,
H_i)} \text{ or }
\frac{(f_{3},H_{j})}{(f_{3},H_{i})}=\alpha _{ij}\frac{(f_{2},H_{j})}{%
(f_{2},H_{i})}  \tag{19}
\end{equation*}
We now prove that $\alpha _{ij}=1$ for all $1\le i\ne j\le 3n-1$. 
Indeed, if
there exists $\alpha _{i_{0}j_{0}}\ne 1$, without loss of generality, 
we may
assume that $\displaystyle \frac{(f_{2},H_{j_{0}})}{(f_{2},H_{i_{0}})}%
=\alpha _{i_{0}j_{0}}\displaystyle 
\frac{(f_{1},H_{j_{0}})}{(f_{1},H_{i_{0}})%
}$. On the other hand $f_{1}=f_{2}$ on $\Omega 
:=\bigcup\limits_{j=1}^{q}%
\big\{z:v_{(f_{1},H_{j})}^{k)}(z)>0\big\}$. Hence, $%
(f_{1},H_{j_{0}})=(f_{2},H_{j_{0}})=0$ on $\Omega \setminus
f_{1}^{-1}(H_{i_{0}})$. So we have
\begin{equation*}
\sum_{j=1,j\neq i_{0}}^{q}\overline{N}_{f_{1}}^{k)}(r,H_{j})\leq 
N\Big(r,v_{%
\frac{(f_{1},H_{j_{0}})}{(f_{1},H_{i_{0}})}}\Big)+\Big(\overline{N}\big(%
r,v_{(f_{1},H_{i_{0}})}\big)-\overline{N}^{k)}\big(r,v_{(f_{1},H_{i_{0}})}%
\big)\Big). 
\end{equation*}
Thus, by the First and the Second Main Theorem, we have
\begin{eqnarray*}
(q-n-2)T_{f_{1}}(r) &\le &\sum_{j=1,j\ne
i_{0}}^{q}N_{n,f_{1}}(r,H_{j})+o(T_{f_{1}}(r)) \\
\hskip3cm &\le &n\sum_{j=1,j\ne
i_{0}}^{q}N_{1,f_{1}}(r,H_{j})+o(T_{f_{1}}(r)) \\
\end{eqnarray*}

\begin{eqnarray*}
&\le &\frac{nk}{k+1}\sum_{j=1,j\ne i_{0}}^{q}\overline{N}%
_{f_{1}}^{k)}(r,H_{j})+\frac{n}{k+1}\sum_{j=1,j\ne
i_{0}}^{q}N_{f_{1}}(r,H_{j})+o(T_{f_{1}}(r)) \\
&\le &\frac{nk}{k+1}N\Big( 
r,v_{\frac{(f_{1},H_{j_{0}})}{(f_{1},H_{i_{0}})}}%
\Big)+\frac{nk}{k+1}\Big( \overline{N}\big( 
r,v_{(f_{1},H_{i_{0}})}\big) -%
\overline{N}^{k)}\big( r,v_{(f_{1},H_{i_{0}})}\big)\Big) \\
&&\hskip3cm+\frac{(q-1)n}{k+1}T_{f_{1}}(r)+o(T_{f_{1}}(r))
\end{eqnarray*}

\begin{align*}
& \le 
\frac{nk}{k+1}T_{\frac{(f_{1},H_{j_{0}})}{(f_{1},H_{i_{0}})}}(r)+\frac{%
nk}{(k+1)^{2}}N_{f_{1}}(r,H_{i_{0}})+\frac{(q-1)n}{k+1}T_{f_{1}}(r)+o\big(%
T_{f_{1}}(r)\big) \\
& \le \Big(\frac{nk}{k+1}+\frac{nk}{(k+1)^{2}}+\frac{(q-1)n}{k+1}\Big) 
T_{f_{1}}(r)+o\big(T_{f_{1}}(r)\big)
\end{align*}
Thus, we get $(q-n-2)\le 
\displaystyle\frac{nk}{k+1}+\displaystyle\frac{nk}{(k+1)^{2}%
}+\displaystyle\frac{(q-1)n}{k+1}\leq n+\frac{qn}{k}$. This contradicts 
any of the following cases:

\noindent i)\ \ $2\leq n\leq 3,\;\;q=3n+1$ and $k\geq 23n$,

\noindent ii)\ \ \thinspace $4\leq n\leq 6,\;\;q=3n$ and $k\geq 
\frac{(6n-1)n%
}{n-3}$,

\noindent iii) $n\geq 7,\;\;q=3n-1$ and $k\geq \frac{(6n-4)n}{n-6}$.\\
\noindent Thus $\alpha_{ij}=1$ for all $1 \le i \ne j \le 3n-1.$

By (19), for $i=3n-1,j\in \{1,\ldots ,3n-2\}$, without loss of 
generality, we
may asssume that 
\begin{equation*}
\frac{(f_{1},H_{j})}{\big(f_{1},H_{3n-1}\big)}=\frac{%
(f_{2},H_{j})}{\big(f_{2},H_{3n-1}\big)}, \;j=1,\ldots ,n: \tag{20}
\end{equation*}
For $1\leq s<v\leq 3$, denote by $L_{sv}$ the set of all $j\in 
\{1,...,3n-2\}$
such that 
$\frac{(f_{s},H_{j})}{(f_{s},H_{3n-1})}=\frac{(f_{v},H_{j})}{%
(f_{v},H_{3n-1})}$ .
By ($19$)  we have  $L_{12}\cup L_{23}\cup L_{13}=\{1,...,3n-2\}.$
So by Dirichlet we have that one of the three sets contains at least 
$n$ different indices,
which are, without loss of generality, $j=1,...,n$, which proves (20).

We choose homogeneous coordinates $(\omega _{0}:\cdots :\omega _{n})$ 
on $\mathbb{C}%
P^{n}$ with $H_{j}=\{\omega _{j}=0\}$ $(1\le j\le n)$, 
$H_{3n-1}=\{\omega
_{0}=0\}$ and take reduced representations: $f_{1}=(f_{1_{0}}:\cdots
:f_{1_{n}})$, $f_{2}=(f_{2_{0}}:\cdots :f_{2_{n}}).$ Then by (20) we 
have
\begin{equation*}
\begin{cases}\displaystyle\frac{f_{1_{j}}}{f_{1_{0}}}=\displaystyle\frac{f_{2_{j}}}{%
f_{2_{0}}}\\
(j=1,\ldots ,n)\end{cases}\;\;\Rightarrow 
\;\;\frac{f_{1_{0}}}{f_{2_{0}}}=\cdots =\frac{f_{1_{n}}}{f_{2_{n}}}\;\Rightarrow f_{1}\equiv f_{2}. 
\end{equation*}
This is a contradiction. Thus, for any case we have that 
$f_{1},f_{2},f_{3}$ 
can not be distinct. Hence, the Proof of Theorem 1 is complete. \hfill  
$\square $\\

\noindent \textbf{Proof of Theorem 2.}
Assume that $\#F_{k}(\{H_{j}\}_{j=1}^{q},f,1)\geq 3.$ Take arbitrarily 
three distinct mappings $f_{1},f_{2},f_{3}\in $\ 
$F_{k}(\{H_{j}\}_{j=1}^{q},f,1). $
We have to prove that $f_{s}\times f_{v}:\mathbb{C}^{m}\longrightarrow 
\mathbb{C}P^{n}\times \mathbb{C}P^{n}$ is linearly degenerate for all 
$1\leq
s<v\leq 3.$

Denote by $Q$ the set which contains all indices $j\in \{1,...,q\}$ 
satisfing 
 $\Phi ^{l}\big(F_{1c}^{j},F_{2c}^{j},F_{3c}^{j}\big)\not\equiv 0$  for 
some $c\in 
\mathcal{C} .$ We distinguish between the two cases $n$ odd and $n$ 
even:\\

\noindent \textbf{Case 1.} If $n$ is odd, then $q=\frac{5(n+1)}{2}.$\\
We now pove that: $Q=$ \ $\emptyset $\ . \hfill (21)

Indeed, otherwise there exist $j_{0}\in Q$ . Then by Lemma 2 (with 
$A=\emptyset
, p=1$) we have
\begin{equation*}
\ 2\sum_{j=1,j\ne j_{0}}^{q}\overline{N}_{f_{i}}^{k)}(r,H_{j})\leq \ \ 
\frac{%
k+2}{k+1}T(r)+o(T(r)),\ i=1,2,3\ . 
\end{equation*}
 (note that $N_{0,f_{i}}^{k)}(r,H_{j_{0}})=0$).
On the other hand, by Lemma 3 we have
\begin{equation*}
2\sum_{j=1,j\ne 
j_{0}}^{q}\overline{N}_{f_{i}}^{k)}(r,H_{j})+o(T_{f_{i}}(r))%
\geq \frac{2[(q-n-2)(k+1)-(q-1)n]}{nk}T_{f_{i}}(r)\ ,\ i=1,2,3\ . 
\end{equation*}
 Hence, we have
$$ \left( (2q-2n-4)(k+1)-2(q-1)n\right) T_{f_{i}}(r)  \leq  
\frac{(k+2)nk}{
k+1}T(r)+o(T(r))\ ,
 i=1,2,3,$$ which implies
\begin{align*}
\Big( (2q-2n-4)(k+1)-2(q-1)n\Big) T(r) & \leq & \frac{
3(k+2)nk}{k+1}T(r)+o(T(r)) \\ 
& \leq & 3n(k+1)T(r)+o(T(r)). 
\end{align*}
Hence, we obtain
$$(2q-2n-4)(k+1)-2(q-1)n\leq  3n(k+1)$$
implying
$$  k+1\leq (5n+3)n. $$   
This is a contradiction. Thus, we get ($21$).\\

\noindent \textbf{Case 2. }If $n$ is even, then $q=\frac{5n+4}{2}$ .\\
We now prove that $\#Q\leq 1$. \hfill (22)\\
Indeed, suppose that this assertion does not hold, then there exist two
distinct indices $j_{0},j_{1}\in Q$ .
By Lemma 2 (with $A=\emptyset ,p=1$) we have
\begin{equation*}
\ \ 2\sum_{j=1,j\ne j_{0}}^{q}\overline{N}_{f_{i}}^{k)}(r,H_{j})\leq 
\frac{%
k+2}{k+1}T(r)+o(T(r))\ \ \ \ \ \ ,\ i=1,2,3, 
\end{equation*}
which implies that, for i=1,2,3 
\begin{eqnarray*}
\ 2\sum_{j=1,j\ne j_{0}}^{q}\left( \ 
\overline{N}_{f_{i}}^{k)}(r,H_{j})-\ 
\frac{1}{n}N_{n,f_{i}}^{k)}(r,H_{j})\right) &\leq &\frac{k+2}{k+1}%
T(r)+o(T(r))\  \\
-\frac{2}{n}\sum_{j=1,j\ne i_{0}}^{q}N_{n,f_{i}}^{k)}(r,H_{j})\; ,\;i 
&=&1,2,3.
\end{eqnarray*}
Hence, we get
\begin{eqnarray*}
2\sum_{j=1,j\ne j_{0}}^{q}\sum_{i=1}^{3}\left( \overline{N}%
_{f_{i}}^{k)}(r,H_{j})-\ \frac{1}{n}N_{n,f_{i}}^{k)}(r,H_{j})\right) 
&\leq &%
\frac{3(k+2)}{k+1}T(r)+o(T(r))\ \  \\
&&-\frac{2}{n}\sum_{j=1,j\ne
j_{0}}^{q}\sum_{i=1}^{3}N_{n,f_{i}}^{k)}(r,H_{j}), \text{(23)}
\end{eqnarray*}
By Lemma 3 (with $q=\frac{5n+4}{2})$, we have 
\begin{equation*}
2\sum_{j=1,j\ne j_{0}}^{q}N_{n,f_{i}}^{k)}(r,H_{j})+o(T_{f_{i}}(r))\geq 
\frac{3n(k+1)-(5n+2)n}{k}T_{f_{i}}(r)\ ,\ i=1,2,3. 
\end{equation*}
Hence, we have
\begin{equation*}
\frac{2}{n}\sum_{j=1,j\ne
j_{0}}^{q}\sum_{i=1}^{3}N_{n,f_{i}}^{k)}(r,H_{j})+o(T(r))\geq \frac{%
3n(k+1)-(5n+2)n}{nk}T(r) \tag{24}
\end{equation*}
By ($23$) and ($24$) we have
\begin{eqnarray*}
2\sum_{j=1,j\ne j_{0}}^{q}\sum_{i=1}^{3}\left( \overline{N}
_{f_{i}}^{k)}(r,H_{j})-\ \frac{1}{n}N_{n,f_{i}}^{k)}(r,H_{j})\right) 
&\leq &%
\frac{(5n+2)n(k+1)-3n}{nk(k+1)}T(r)+o(T(r)) \\
&\leq &\frac{5n+2}{k}T(r)+o(T(r)).
\end{eqnarray*}
On the other hand, we obtain 
\begin{equation*}
\overline{N}_{f_{i}}^{k)}(r,H_{j})-\ 
\frac{1}{n}N_{n,f_{i}}^{k)}(r,H_{j})%
\geq 0\text{ for all }i\in \{1,2,3\},j\in \{1,...,q\}. 
\end{equation*}
Hence, we get
\begin{equation*}
\sum_{i=1}^{3}\left( \overline{N}_{f_{i}}^{k)}(r,H_{j})-\ \frac{1}{n}%
N_{n,f_{i}}^{k)}(r,H_{j})\right) \text{ }\leq 
\frac{5n+2}{k}T(r)+o(T(r))%
\text{ , }j\in \{1,...,q\}\setminus \{j_{0}\}. 
\end{equation*}
In particular, we get
\begin{equation*}
\sum_{i=1}^{3}\left( \text{ }\overline{N}_{f_{i}}^{k)}(r,H_{j_{1}})-\ 
\frac{1%
}{n}N_{n,f_{i}}^{k)}(r,H_{j_{1}})\right) \text{ }\leq \frac{5n+2}{k}%
T(r)+o(T(r)) \tag{25}
\end{equation*}

Set $A_{i}:=\{z\in \mathbb{C}^{m}:v_{(f_{i},H_{j_{1}})}(z)=1\}$ for 
$i=1,2,3$ . For each $i\in \{1,2,3\},$ we have $\overline{A}_{i}\setminus 
A_{i}\subseteq 
$ ${\rm sing}f_{i}^{-1}(H_{j_{1}}).$ Indeed, otherwise there existed 
$a\in \left( \overline{A}_{i}\setminus A_{i}\right) \cap {\rm 
reg}f_{i}^{-1}(H_{j_{1}})$ . Then
 $p_{0}:=v_{(f_{i},H_{j_{1}})}(a)\geq 2$. Since $a$ is a regular point 
of $%
f_{i}^{-1}(H_{j_{1}})$ we can choose nonzero holomorphic functions $h$ 
and $u$
on a neighborhood $U$ of $a$ such that d$h$ and $u$ have no zeroes and 
$%
(f_{i},H_{j_{1}})\equiv h^{p_{0}}u$ on $U$. Since $a\in $ 
$\overline{A}_{i}$ there 
exists $b\in A_{i}\cap U$ . 
Then, we get $1=v_{(f_{i},H_{j_1})}(b)=v_{h^{p_{0}}u}(b)=p_{0}\geq 2.$ 
This is a contradiction.
Thus, we get that  $\overline{A}_{i}\setminus A_{i}\subseteq $ $
{\rm sing}f_{i}^{-1}(H_{j_{1}}). $ 

Set $B:=A_{1}\cup A_{2}\cup A_{3}$. Then $%
\overline{B}\setminus B\subseteq
\bigcup\limits_{i=1}^{3}{\rm sing}f_{i}^{-1}(H_{j_{1}})$. 
This means that $\overline{B}\setminus B$ is included in an analytic 
set of codimension 
$\geq 2$. So we have
\begin{equation*}
(n-1)\overset{\_\_}{N}(r,\overline{B})\leq \sum_{i=1}^{3}\left( n\text{ 
}%
\overline{N}_{f_{i}}^{k)}(r,H_{j_{1}})-\
N_{n,f_{i}}^{k)}(r,H_{j_{1}})\right) . 
\end{equation*}
By (25) we have
\begin{equation*}
\overset{\_\_}{N}(r,\overline{B})\leq 
\frac{(5n+2)n}{(n-1)k}T(r)+o(T(r)),
\end{equation*}
where we note that $n\geq 2$ , since $n$ is even.
It is clear that \\
\textit{min}$\big\{v_{(f_{1},H_{j_{1}})}^{k)},2\big\}=\textit{min}\big\{%
v_{(f_{2},H_{j_{1}})}^{k)},2\big\}=\textit{min}\big\{v_{(f_{3},H_{j_{1}})}^{k)},2%
\big\}$ on $\mathbb{C}^{m}\setminus \overline{B}(\subseteq 
\mathbb{C}^{m}\diagdown B).$
By Lemma 2 (with $A=\overline{B}$ , $p=2$) we have
\begin{eqnarray*}
2\sum_{j=1,j\ne 
j_{1}}^{q}\overline{N}_{f_{i}}^{k)}(r,H_{j})+\overline{N}%
_{f_{i}}^{k)}(r,H_{j_{1}}) &\leq &\frac{k+2}{k+1}T(r)+4\overline{N}(r,%
\overline{B})+o(T(r)) \\
&\leq &\left( \frac{k+2}{k+1}+\frac{4(5n+2)n}{(n-1)k}\right) 
T(r)+o(T(r))%
\text{ , (26)}
\end{eqnarray*}
(note that $j_{1}\in Q$ ).
By Lemma 3 we have
\begin{eqnarray*}
\sum_{j=1,j\ne 
j_{1}}^{q}\overline{N}_{f_{i}}^{k)}(r,H_{j})+o(T_{f_{i}}(r))
&\geq &\frac{(q-n-2)(k+1)-(q-1)n}{nk}T_{f_{i}}(r)\text{ },\text{and} \\
\sum_{j=1}^{q}\overline{N}_{f_{i}}^{k)}(r,H_{j})+o(T_{f_{i}}(r)) &\geq 
&%
\frac{(q-n-1)(k+1)-qn}{nk}T_{f_{i}}(r)\text{ .}
\end{eqnarray*}
Consequently, we obtain  
\begin{align*}
2\sum_{j=1,j\ne 
j_{1}}^{q}\overline{N}_{f_{i}}^{k)}(r,H_{j})+\overline{N}%
_{f_{i}}^{k)}(r,H_{j_{1}})+o(T_{f_{i}}(r))\geq 
\frac{(2q-2n-3)(k+1)-(2q-1)n}{%
nk}T_{f_{i}}(r)\\  (27) 
\end{align*}
By ($26$) and ($27$) we have
$$ \frac{(2q-2n-3)(k+1)-(2q-1)n}{nk}T_{f_{i}}(r)  \leq  \left( 
\frac{k+2}{%
k+1}+\frac{4(5n+2)n}{(n-1)k}\right) T(r)+o(T(r)) ,$$
which implies
\begin{align*}
\left( (3n+1)(k+1)-(5n+3)n\right) T(r) & \leq & \left( \frac{%
3nk(k+2)}{k+1}+\frac{12(5n+2)n^{2}}{(n-1)}\right) T(r)+o(T(r)) \\ 
  & \leq & (3n(k+1)+\frac{12(5n+2)n^{2}}{(n-1)})T(r)+o(T(r)),
\end{align*}
and, hence,
$$ k+1\leq (5n+3)n+\frac{12(5n+2)n^{2}}{(n-1)}.$$
This contradicts $k\geq (65n+171)n$ , $n\geq 2.$ Hence, we have  
$\#Q\leq 1.$ 
So  we get (22).\\

By ($21$) and ($22$) we have $\#(\{1,...,q\}\setminus Q)\geq q-1$ . 
Without loss of generality we may assume that $1,...,q-1\notin Q$ .
For any $j\in \{1,\ldots ,q-1\}$ we have

$\Phi^l\big(F^j_{1c}, F^j_{2c}, F^j_{3c}\big)\equiv 0$ for all $c \in 
\mathcal{C} $, $l
\in \{1,..., m\}$.

\noindent On the other hand, $%
\mathcal{C}$ is dense in $\mathbb{C}^{n+1}$ . Hence, we get that $\Phi 
^{l}\big(
F_{1c}^{j},F_{2c}^{j},F_{3c}^{j}\big )\equiv 0$ for all $c\in 
\mathbb{C}%
^{n+1}\setminus \{0\}$, $l\in \{1,...,m\},$ $j\in \{1,\ldots ,q-1\}$ . 
In particular (for $H_{c}=H_{i})$, we get
\begin{equation*}
\Phi ^{l}\left( 
\frac{(f_{1},H_{j})}{(f_{1},H_{i})},\;\;\frac{(f_{2},H_{j})}{%
(f_{2},H_{i})},\;\;\frac{(f_{3},H_{j})}{(f_{3},H_{i})}\right) \equiv 0 
\end{equation*}
for all $1\le i\ne j\le q-1,\;\;\;l\in \{1,\ldots ,m\}$ .

For each $1\le i\ne j\le q-1,$ by Theorem 2.2, there exists a constant 
$
\alpha _{ij}$ such that
\begin{equation*}
\frac{(f_2, H_j)} {(f_2, H_i)} =\alpha_{ij} \frac{(f_1, H_j)} {(f_1, 
H_i)}
\;{\rm or }\; \frac{(f_3, H_j)} {(f_3, H_i)}=\alpha_{ij} \frac{(f_1, 
H_j)} {(f_1,H_i)} \;{\rm or }\;\frac{(f_{3},H_{j})}{(f_{3},H_{i})}=\alpha 
_{ij}\frac{(f_{2},H_{j})}{(f_{2},H_{i})}\text{ .} 
\end{equation*}
We now prove that \\
$\alpha _{ij}=1$ for all $1\le i\ne j\le q-1$. \hfill (28)\\
Indeed, if there exists $\alpha _{i_{0}j_{0}}\ne 1$, without loss of 
generality, we may assume that $\displaystyle \frac{(f_{2},H_{j_{0}})}{%
(f_{2},H_{i_{0}})}=\alpha _{i_{0}j_{0}}\displaystyle 
\frac{(f_{1},H_{j_{0}})%
}{(f_{1},H_{i_{0}})}$. On the other hand, we have $f_{1}=f_{2}$ on $%
D:=\bigcup\limits_{j=1}^{q}\big\{z:v_{(f_{1},H_{j})}^{k)}>0\big\}$. 
Hence, we get $%
(f_{1},H_{j_{0}})=(f_{2},H_{j_{0}})=0$ on $D\setminus
f_{1}^{-1}(H_{_{i_{0}}})$. So we have
\begin{equation*}
\sum_{j=1,j\neq i_{0}}^{q}\overline{N}_{f_{1}}^{k)}(r,H_{j})\leq 
N\Big(r,v_{%
\frac{(f_{1},H_{j_{0}})}{(f_{1},H_{i_{0}})}}\Big)+\Big(\overline{N}\big(%
r,v_{(f_{1},H_{i_{0}})}\big)-\overline{N}^{k)}\big(r,v_{(f_{1},H_{i_{0}})}%
\big)\Big). 
\end{equation*}
Thus, by the First and the Second Main Theorem, we have 
\begin{eqnarray*}
(q-n-2)T_{f_{1}}(r) &\le &\sum_{j=1,j\ne
i_{0}}^{q}N_{n,f_{1}}(r,H_{j})+o(T_{f_{1}}(r)) \\
\hskip3cm &\le &n\sum_{j=1,j\ne
i_{0}}^{q}N_{1,f_{1}}(r,H_{j})+o(T_{f_{1}}(r))
\end{eqnarray*}
\begin{eqnarray*}
&\le &\frac{nk}{k+1}\sum_{j=1,j\ne i_{0}}^{q}\overline{N}%
_{f_{1}}^{k)}(r,H_{j})+\frac{n}{k+1}\sum_{j=1,j\ne
i_{0}}^{q}N_{f_{1}}(r,H_{j})+o(T_{f_{1}}(r)) \\
&\le &\frac{nk}{k+1}N\Big( 
r,v_{\frac{(f_{1},H_{j_{0}})}{(f_{1},H_{i_{0}})}}%
\Big)+\frac{nk}{k+1}\Big( \overline{N}\big( 
r,v_{(f_{1},H_{i_{0}})}\big) -%
\overline{N}^{k)}\big( r,v_{(f_{1},H_{i_{0}})}\big)\Big) \\
&&\hskip3cm+\frac{(q-1)n}{k+1}T_{f_{1}}(r)+o(T_{f_{1}}(r))
\end{eqnarray*}
\begin{eqnarray*}
&\le 
&\frac{nk}{k+1}T_{\frac{(f_{1},H_{j_{0}})}{(f_{1},H_{i_{0}})}}(r)+\frac{%
nk}{(k+1)^{2}}N_{f_{1}}(r,H_{i_{0}})+\frac{(q-1)n}{k+1}T_{f_{1}}(r)+o\big(%
T_{f_{1}}(r)\big) \\
&\le &\Big(\frac{nk}{k+1}+\frac{nk}{(k+1)^{2}}+\frac{(q-1)n}{k+1}\Big) 
T_{f_{1}}(r)+o\big(T_{f_{1}}(r)\big)
\end{eqnarray*}
Thus, we have  $(q-n-2)\le 
\displaystyle\frac{nk}{k+1}+\displaystyle\frac{nk}{(k+1)^{2}}+\displaystyle\frac{(q-1)n}{k+1}\leq n+\frac{nq}{k}$. 
\\
This contradicts $q=\left[ \frac{5(n+1)}{2}\right] , k\geq (65n+171)n.$
Thus, we get that $\alpha _{ij}=1$ for all $1\le i\ne j\le q-1.$\\

For $1\leq s<v\leq 3$, denote by $L_{sv}$ the set of all $j\in 
\{1,...,q-2\}$
such that $\frac{(f_{s},H_{j})}{(f_{s},H_{q-1})}=\frac{(f_{v},H_{j})}{%
(f_{v},H_{q-1})}$ .
By ($28$) , we have that $L_{12}\cup L_{23}\cup L_{13}=\{1,...,q-2\}.$

 If there exists some $L_{sv}=\emptyset $, without loss of generality, 
we may
 assume that $L_{13}=\emptyset .$ Then $L_{12}\cup 
L_{23}=\{1,...,q-2\}.$
Since $q=\left[ \frac{5(n+1)}{2}\right] $ we have that $\#L_{12}\geq n$ 
or $%
\#L_{23}\geq n.$ We may assume that $\#L_{12}\geq n$ , and furthermore 
$%
1,...,n\in L_{12}$ . Then 
$\frac{(f_{1},H_{j})}{(f_{1},H_{q-1})}=\frac{%
(f_{2},H_{j})}{(f_{2},H_{q-1})}$  for all $j\in \{1,...,n\}$, so 
$f_{1}\equiv f_{2}$ (as in the proof of Theorem 1). This is a contradiction.

Thus, we have $L_{sv}\neq \emptyset $ for all $1\leq s<v\leq 3.$ Then 
for any $1\leq s<v\leq 3,$ there exists $j\in \{1,...,q-2\}$ such that 
$\frac{(f_{s},H_{j})}{%
(f_{s},H_{q-1})}=\frac{(f_{v},H_{j})}{(f_{v},H_{q-1})}.$
Hence, we finally get that   $f_{s}\times 
f_{v}:\mathbb{C}^{m}\longrightarrow \mathbb{C}P^{n}\times 
\mathbb{C}P^{n}$ is linearly degenerate. We thus have completed the 
proof of Theorem 2. \hfill  $\square $\\

 \noindent Gerd Dethloff  \\
Universit\'{e} de Bretagne Occidentale \\
    UFR Sciences et Techniques \\
D\'{e}partement de Math\'{e}matiques \\
6, avenue Le Gorgeu, BP 452 \\
   29275 Brest Cedex, France \\
e-mail: gerd.dethloff@univ-brest.fr\\

\noindent Tran Van Tan\\
Department of Mathematics\\
  Hanoi University of Education\\
 Cau Giay, Hanoi, Vietnam\\
e-mail: tranvantanhn@yahoo.com

\end{document}